\documentclass[12pt,leqno]{article}
\usepackage{amssymb,eufrak}
\begin{document}
\title{On some $p$-adic power series attached to the arithmetic of
$\mathbb Q(\zeta_p )$}

\author{
Bruno Angl\`  es \\
Universit\'e de Caen, \\
Laboratoire Nicolas Oresme, CNRS UMR 6139, \\
Campus II, Boulevard Mar\'echal Juin, \\
BP 5186, 14032 Caen Cedex, France.\\
E-mail: angles@math.unicaen.fr}
\date{October 11  2005}

\maketitle
Let $p$ be a prime number, $p\geq 5.$ Let $\theta$ be an even and non-trivial character of $\rm{Gal}(\mathbb Q(\mu_p )/\mathbb Q ).$
Let $f(T,\theta )\in \mathbb Z_p[[T]]$ be the Iwasawa power series associated to the $p$-adic L-function $L_p(s,\theta )$ (see
\cite{WAS}), i.e.:
$$\forall n\geq 1, n\equiv 0\pmod{p-1},\, f((1+p)^{1-n}-1,\theta )=L(1-n, \theta ),$$
where $L(s, \theta )$ is the usual Dirichlet L-series. In 1979, in their celebrated article \cite{FER}, Ferrero and Washington have
proved that (see also \cite{ADL}):
$$f(T, \theta )\not \equiv 0\pmod{p}.$$
Thus we can write:
$$f(T, \theta )\equiv T^{\lambda (\theta )} \bar {U}(T)\pmod{p},$$
where $\bar{U}(T)\in \mathbb F_p[[T]]^*.$  The Iwasawa lambda-invariant $\lambda (\theta )$ is not well-understood. By a heuristic
argument due to Ferrero and Washington, one could expect that for sufficiently large $p$ (see \cite{LAN}):
$$\sum_{\theta \, \rm{even}\, ,\theta \not = 1}\lambda (\theta )\leq \frac{{\rm Log}(p)}{{\rm Log}({\rm Log}(p))}.$$
Furthermore, if $p< 4 000 000,$ we have $\lambda (\theta )\leq 1 ,$ and one would reasonably expect (see \cite{MET} ):
$$\lambda(\theta )< p.$$
The only known result is due to Ferrero and Washington (\cite{FER}):  for sufficiently large $p$ we have:
 $$\lambda (\theta )\leq
p^{\rm{Log}(p)^{4\varphi (p-1)^4}}.$$ 
Now observe that:
$$f'(T, \theta )\equiv \lambda (\theta )T^{\lambda (\theta )-1} \bar{U}(T) +T^{\lambda (\theta )}\bar {U}'(T)\pmod{p}.$$
Thus if $\lambda (\theta )\geq 1$ and $\lambda (\theta ) \not \equiv 0 \pmod{p},$ we have that $f'(T, \theta )\not \equiv 0
\pmod{p}.$\par
The aim of this paper is to prove that for all $\theta $ even, $\theta \not =1,$ we have (Corollary \ref{Corollary1}):
$$f'(T, \theta )\not \equiv 0 \pmod{p}.$$
This fact comes from  properties of some power series that are connected to polynomials introduced by Mirimanoff at the beginning of the
XXth century.\par

\section{Notations}\par
${}$\par
Let $p$ be a prime number, $p\geq 5.$ Let $\overline{\mathbb Q_p}$ be an
algebraic closure of $\mathbb Q_p .$ All the extensions of $\mathbb
Q_p$ considered in this paper are contained in $\overline{\mathbb
Q_p}.$ Let $v_p$ be the $p$-adic valuation on $\overline{\mathbb Q_p}$ such that $v_p(p)=1.$ We denote the
Iwasawa $p$-adic logarithm on $\overline{\mathbb Q_p}$ by ${\rm Log}_p.$\par
${}$\par
If $A$ is a commutative ring, we denote the set of invertible
elements of $A$ by $A^*.$\par
${}$\par
 For every integer $d,$ $d\geq 1 ,$ set $\mu_d=\{ z\in
\overline{\mathbb Q_p} \mid z^d=1\}.$ If $\rho \in \cup_{d\geq 1}\mu_d ,$
we denote the order of $\rho $ by $o(\rho ).$ Set
$\mu_{p^{\infty}}=\cup_{n\geq 0}\mu_{p^{n+1}}.$ For all $n\geq 0,$ let
$\zeta_{p^{n+1}}\in \mu_{p^{n+1}}$ such that: 
\noindent $\zeta_p\not = 1$ and $\forall n\geq 0,$
$\zeta_{p^{n+2}}^p=\zeta_{p^{n+1}}.$\par
${}$\par
Let:\par
\noindent - $K_n=\mathbb Q_p(\mu_{p^{n+1}} ),$\par
\noindent - $O_n=\mathbb Z_p[\zeta_{p^{n+1}}],$\par
\noindent - $\pi_n=\zeta_{p^{n+1}}-1,$\par
\noindent - $U_n=1+\pi_nO_n,$\par
\noindent - $\Gamma_n={\rm{Gal}}(K_n/K_0),$\par
\noindent - $\Delta ={\rm{Gal}}(K_0/\mathbb Q_p),$\par
\noindent - $K_{\infty}=\mathbb Q_p(\mu_{p^{\infty }}),$\par
\noindent - $\Gamma ={\rm{Gal}}(K_{\infty }/K_0).$\par
\noindent Let $\gamma_0 \in \Gamma $ be such that $\forall \zeta \in
\mu_{p^{\infty }},$ $\zeta^{\gamma_0}=\zeta^{1+p}.$\par
\noindent For $a \in \mathbb Z_p^*,$ we write:
$$a=\omega (a) <a>,$$ where $\omega $ is the Teichm\"uller character
(i.e. $\omega (a)=\rm{lim}_n a^{p^n}\, \in \mu_{p-1}$) and
$<a>\equiv 1 \pmod{p}.$\par
\noindent Let $n\geq 0.$ For $a\in \mathbb Z_p^*,$ let $\sigma_n (a)\in
{\rm{Gal}}(K_n/\mathbb Q_p)$ be such that:
$$\zeta_{p^{n+1}}^{\sigma_n(a)}=\zeta_{p^{n+1}}^a.$$
We have $\sigma_n(a)=\delta (a) \gamma_n(a),$ where $\delta (a) \in
\Delta $ and $\gamma_n(a) \in \Gamma _n.$ Observe that:
$$\zeta_{p^{n+1}}^{\gamma_n(a)}=\zeta_{p^{n+1}}^{<a>}.$$
For $a\in \mathbb Z_p^*,$ let $\sigma (a)\in \Gamma \times \Delta $ be
such that:
$$\forall \zeta \in \mu_{p^{\infty }},\, \zeta^{\sigma (a)} =\zeta^a .$$
Note that $\Delta =\{ \sigma (\omega (a))\mid  a= 1, \cdots , p-1 \}.$ 
For $\theta \in \widehat{\Delta }=\rm{Hom}(\Delta , \mu_{p-1}),$ we set:
$$e_{\theta }=\frac{1}{p-1}\sum_{\delta \in \Delta}\theta ^{-1}(\delta
)\, \delta \, \in \mathbb Z_p[\Delta ].$$
${}$\par
Let $T$ be an indeterminate over $\overline{\mathbb Q_p }.$  Let $L$ be a
finite extension of $\mathbb Q_p .$  Let $O_L$ be the valuation ring of
$L.$ The restriction maps $Res_{n+1,n}:\Gamma_{n+1}\rightarrow \Gamma_n
$ induce maps $Res_{n+1,n}: O_L[\Gamma_{n+1}]\rightarrow O_L[\Gamma _n].$
Thus we can take the inverse limit of the group rings $O_L[\Gamma_n]$
with respect to these maps, and by \cite{WAS}, Theorem 7.1,we have:
$$O_L[[T]]\simeq \lim_{\leftarrow}O_L[\Gamma_n],$$ 
where $T$ corresponds to
$\gamma_0-1.$ We set $\Lambda_L=O_L[[T]]$ and  $\Lambda =\Lambda_{\mathbb Q_p}.$ For all $n\geq 0,$
we set $\omega_n(T)=(1+T)^{p^n}-1.$ Recall that:
$$\forall n\geq 0,\, \frac{\Lambda_L}{\omega_n(T)\Lambda_L}\simeq
O_L[\Gamma_n].$$
${}$\par

Note that forl all $n\geq 0,$ $N_{K_{n+1}/K_n}(U_{n+1})\subset U_n,$
where $N_{K_{n+1}/K_n}$ is the norm map from $K_{n+1}$ to $K_n.$ We
denote the inverse limit of the principal units $U_n$ with respect to the
norm maps by $U_{\infty }.$ Note that $U_{\infty }$ is a $\Lambda [\Delta
]$-module.\par
\section{Logarithmic derivatives}\par
${}$\par
Let $u=(u_n)_{n\geq 0}$ be an element in $U_{\infty },$ recall that:
$$\forall n\geq 0,\, N_{K_{n+1}/K_n}(u_{n+1})=u_n.$$
There exists an unique element $f_u(T) \in \Lambda $ such that (\cite{WAS}, Theorem 13.38, see also
\cite{COL}):
$$\forall n\geq 0,\, f_u(\pi_n)=u_n.$$
Furthermore, if $x=\sum b_a \sigma (a)\in \mathbb Z_p[\Gamma \times \Delta ]$  and if $u\in U_{\infty },$ then
we have (\cite{WAS}, Lemma 13.48):
$$f_{u^x}(T)=\prod f_u((1+T)^{a} -1)^{b_a}.$$
For $u\in U_{\infty }$ and for $n\geq 0,$ we set:
$$D_n(u)=\zeta_{p^{n+1}}\frac {f_u'(\pi_n)}{f_u(\pi_n)} \, \in O_n.$$
We call the map $D_n: U_{\infty } \rightarrow O_n$ the $n$th logarithmic derivative of Coleman.\par
\newtheorem{Lemma1}{Lemma}[section]
\begin{Lemma1} \label{Lemma1}
${}$\par
\noindent i) Let $u\in U_{\infty },$ $\forall n\geq 0,$ $D_n(u^{\omega_n(T)})\equiv 0\pmod{p^{n+1}}.$\par
\noindent ii) $\forall \theta \in \widehat{\Delta },$ $\forall u\in U_{\infty },$ $\forall n\geq 0,$
$D_n(u^{e_{\theta }})=e_{\theta \omega^{-1}}D_n(u).$\par
\end{Lemma1}
\noindent{\sl Proof}${}$\par
\noindent i) We have:
$$f_{u^{\omega_n(T)}}(T)=\frac{f_u((1+T)^{(1+p)^{p^n}}-1)}{f_u(T)}.$$
Thus:
$$(1+T)\frac{f_{u^{\omega_n(T)}}'(T)}{f_{u^{\omega_n(T)}}(T)}= (1+p)^{p^n}(1+T)^{(1+p)^{p^n}}\frac
{f_u'((1+T)^{(1+p)^{p^n}}-1)}{f_u((1+T)^{(1+p)^{p^n}}-1)}-(1+T)\frac{f_u'(T)}{f_u(T)}.$$
Note that:
$$(1+\pi_n)^{(1+p)^{p^n}}-1=\pi_n.$$
Therefore:
$$D_n(u^{\omega_n(T)})=\zeta_{p^{n+1}}((1+p)^{p^n}-1)\frac{f_u'(\pi_n)}{f_u(\pi_n)}\equiv 0\pmod{p^{n+1}}.$$
ii) Let $a\in \mathbb Z_p^*$ and let $u\in U_{\infty }.$ Observe that:
$$D_n(u^{\sigma (a)})=a\sigma_n(a) (D_n(u)).$$
Let $u\in U_{\infty },$ we get:
$$D_n(u^{e_{\theta }})=\frac{1}{p-1}\sum_{a=1}^{p-1}\theta^{-1}(a) D_n(u^{\sigma (\omega (a))}).$$
We have:
$$D_n(u^{e_{\theta }})=\frac{1}{p-1}\sum_{a=1}^{p-1}\theta^{-1}(a)\, \omega (a) \, \sigma_n(\omega (a)) (D_n
(u)).$$ The Lemma follows. $\diamondsuit$\par
\newtheorem{Proposition1}[Lemma1]{Proposition}
\begin{Proposition1} \label{Proposition1}
Let $\theta \in \widehat {\Delta },$ $\theta \not = 1,\omega.$  The map $D_n$ gives rise to a morphism of
$\mathbb Z_p[\Gamma_n]$-modules:
$$\phi_n^{\theta }: \frac{U_n^{e_{\theta }}}{U_n^{p^{n+1}e_{\theta }}} \rightarrow \frac{e_{\theta
\omega^{-1}}O_n}{p^{n+1}e_{\theta \omega^{-1}}O_n},$$
where $T$ acts on the right via $(1+p)(1+T)-1.$\par
\end{Proposition1}
\noindent{\sl Proof} We have (\cite{WAS}, Theorem 13.54):
$$U_n^{e_{\theta }}\simeq \frac{U_{\infty }^{e_{\theta }}}{U_{\infty }^{\omega_n(T)e_{\theta }}}.$$
Let $v\in U_n^{e_{\theta }},$ then there exists $u=(u_n)_{n\geq 0}\in U_{\infty }^{e_{\theta }}$ such that
$v=u_n.$ We set:
$$\phi_n^{\theta }(v)\equiv D_n(u)\pmod{p^{n+1}}.$$
By Lemma \ref{Lemma1}, we get a morphism of $\mathbb Z_p$-modules:
$$\phi_n^{\theta }: \frac{U_n^{e_{\theta }}}{U_n^{p^{n+1}e_{\theta }}} \rightarrow \frac{e_{\theta
\omega^{-1}}O_n}{p^{n+1}e_{\theta \omega^{-1}}O_n}.$$
Let $v\in U_n^{e_{\theta }},$ and let $u=(u_n)_{n\geq 0}\in U_{\infty }^{e_{\theta }}$ such that $v=u_n.$ We
have:
$$\phi_n^{\theta }(u^T)\equiv D_n (u^{\gamma_0-1})\pmod{p^{n+1}}.$$
But:
$$D_n(u^{\gamma_0-1})\equiv ((1+p)\gamma_0-1)(D_n(u))\pmod{p^{n+1}}.$$
The Proposition follows. $\diamondsuit$\par
\newtheorem{Theorem1}[Lemma1]{Theorem}
\begin{Theorem1}\label{Theorem1}
For all $n\geq 0,$ set $T_n=\sum_{d=0}^{n}\zeta_{p^{d+1}}.$ Let $\theta \in \widehat{\Delta },$ $\theta \not =
1,\omega .$ Then:
$$\forall n\geq 0,\, \rm{Im}(\phi_n^{\theta })=\frac{\mathbb Z_p[\Gamma_n]e_{\theta
\omega^{-1}}T_n}{p^{n+1}e_{\theta \omega^{-1}}O_n}.$$
\end{Theorem1}
\noindent{\sl Proof} Let $\Lambda_n$ be the maximal order of $\mathbb Q_p[\Gamma_n].$ Then (\cite {LEO}, see
also \cite{LET}):
$$O_n=\Lambda_n[\Delta ]\,  T_n.$$ 
Now observe that $p^{n+1} \Lambda_n \subset p\mathbb Z_p[\Gamma _n].$ Thus:
$$ p^{n+1}e_{\theta \omega^{-1}}O_n \subset p\mathbb Z_p[\Gamma _n]e_{\theta \omega^{-1}} \, T_n.$$
Recall that there exists $\alpha \in \mu_{p-1}\setminus \{ 1\} $ such that (\cite{WAS}, Corollary 13.37):
$$\forall n\geq 0,\, U_n^{e_{\theta}} =(\frac{\alpha -\zeta_{p^{n+1}}}{\omega (\alpha -1)})^{\mathbb Z_p[\Gamma
_n] e_{\theta }}.$$
For $\alpha \in \mu_{p-1}\setminus \{ 1\} ,$ set:
$$\rho_{\infty }(\alpha )=(\frac{\alpha -\zeta_{p^{n+1}}}{\omega (\alpha -1)})_{n\geq 0}\in U_{\infty }.$$
We have:
$$D_n(\rho_{\infty }(\alpha ))=\frac{-\zeta_{p^{n+1}}}{\alpha -\zeta_{p^{n+1}}}.$$
Thus  $\rm{Im}(\phi_n^{\theta })$ is the $\mathbb Z_p[\Gamma_n]$-module generated by the $e_{\theta
\omega^{-1}} \frac{-\zeta_{p^{n+1}}}{\alpha -\zeta_{p^{n+1}}},\, \alpha \in\mu_{p-1}\setminus \{ 1\} .$ For
$\alpha \in \mu_{p-1}\setminus \{ 1\} ,$ for all $n\geq 0,$ set:
$$u_n(\alpha )=\sum_{k=1, \, k\not \equiv 0\pmod{p}}^{p^{n+1}-1} \alpha ^k \sigma_n(k) \in \mathbb Z_p[\Gamma
_n\times \Delta ].$$
We have:
$$Res_{n+1,n} (u_{n+1}(\alpha ))=\sum_{k=1, \, k\not \equiv 0\pmod{p}}^{p^{n+2}-1} \alpha ^k \sigma_n(k).$$
Thus :
$$Res_{n+1,n} (u_{n+1}(\alpha ))=\sum_{k=1, \, k\not \equiv 0\pmod{p}}^{p^{n+1}-1}\sigma_n (k)\, (\sum_{\ell
=0}^{p-1} \alpha^{k+p^{n+1}\ell }).$$
Finally, we get:
$$Res_{n+1,n} (u_{n+1}(\alpha ))=u_{n}(\alpha ).$$
Now, we have:
$$\frac{\alpha -1}{\alpha \zeta_{p^{n+1}}-1}-1=\frac{(\alpha \zeta_{p^{n+1}})^{p^{n+1}}-1}{\alpha
\zeta_{p^{n+1}}-1}-1.$$
Thus:
$$\frac{\alpha -1}{\alpha \zeta_{p^{n+1}}-1}-1=\sum_{k=1}^{p^{n+1}-1}\alpha^k \zeta_{p^{n+1}}^k.$$
Therefore:
$$\frac{\alpha -1}{\alpha \zeta_{p^{n+1}}-1}-1=\sum_{k=1, k\not \equiv 0 \pmod{p}}^{p^{n+1}-1}\alpha^k
\sigma_n(k)(\zeta_{p^{n+1}}) +\sum_{k=1}^{p^{n}-1}\alpha^k \zeta_{p^{n}}^k.$$
Finally, we get:
$$\frac{\alpha -1}{\alpha \zeta_{p^{n+1}}-1}-1=u_n(\alpha )(T_n).$$
Now set:
$$u_n (\theta , \alpha )=\sum_{k=1, \, k\not \equiv 0\pmod{p}}^{p^{n+1}-1}\theta  (k) \omega^{-1}(k) \alpha ^k
\gamma_n(k).$$
We have:
$$(u_n (\theta , \alpha ))_{n\geq 0}\in \Lambda .$$
Furthermore, we have:
$$e_{\theta \omega^{-1}}\frac{-\zeta_{p^{n+1}}}{\alpha -\zeta_{p^{n+1}}}= e_{\theta \omega^{-1}}\frac {1}
{\alpha ^{-1}\zeta_{p^{n+1}}-1}.$$
Thus:
$$e_{\theta \omega^{-1}}\frac{-\zeta_{p^{n+1}}}{\alpha -\zeta_{p^{n+1}}}=\frac{1}{\alpha^{-1}-1}u_n(\theta ,
\alpha^{-1})(e_{\theta \omega^{-1}} T_n).$$
Now:
$$\sum_{\alpha \in \mu_{p-1}\setminus \{ 1\} }u_0(\theta , \alpha)=(p-1) \theta (-1) \omega^{-1} (-1)\in
\mathbb Z_p^*.$$
Thus there exists $\alpha \in \mu_{p-1}\setminus \{ 1 \} $ such that for all $n\geq 0,$ $u_n(\theta ,\alpha )
\in \mathbb Z_p[\Gamma_n]^*.$ The theorem follows. $\diamondsuit$\par
We will also need the following Lemma:
\newtheorem{Lemma2}[Lemma1]{Lemma}
\begin{Lemma2} \label{Lemma2} Let $\Lambda_n$ be the maximal order of $\mathbb Q_p[\Gamma_n].$ Under the
isomorphism $\Lambda/\omega_n(T)\Lambda\simeq \mathbb Z_p[\Gamma _n],$ $p^{n+1}\Lambda_n$ corresponds to an
ideal ${\cal U}_n$ of $\Lambda$ such that $\omega_n(T) \Lambda \subset {\cal U}_n \subset (p, \omega_n(T))$
and $\lim_n {\cal U}_n = \{ 0\} .$
\end{Lemma2}
\noindent{\sl Proof} Set:
$$e_0=\frac{1}{p^n} \sum_{\gamma \in \Gamma_n}\gamma ,$$
and for $1\leq d\leq n ,$ set:
$$e_d=\sum_{\chi \in \widehat{\Gamma_n }, \, f_{\chi }= p^{d+1}}e_{\chi},$$
where the sum is over all the characaters of $\Gamma_n$ of conductors $p^{d+1}$ and :
$$e_{\chi }=\frac{1}{p^n}\sum_{\gamma \in \Gamma_n} \chi^{-1}(\gamma )\gamma .$$
We have:
$$p^{n+1} \Lambda_n=\oplus_{d=0}^{n}\mathbb Z_p[\Gamma_n]p^{n+1}e_d\subset p\mathbb Z_p[\Gamma_n].$$
Thus it is clear that $\omega_n(T)\Lambda \subset {\cal U}_n\subset (p, \omega_n(T)).$ Now:
$$p^{n+1}e_0\equiv p(\sum_{\ell =0}^{p^n -1} (1+T)^{\ell } )\pmod{\omega_n(T)}.$$
Thus:
$$p^{n+1}e_0\equiv p\frac{\omega_n(T)}{T} \pmod{\omega_n (T)}.$$
Let $d$ be an integer, $1\leq d \leq n.$ We have:
$$p^{n+1}e_d\equiv p(\sum_{\ell =0}^{p^n-1}(1+T)^{\ell} (\sum_{f_{\chi }=p^{d+1}} \chi(1+p)^{\ell
}))\pmod{\omega_n(T)}.$$
Thus:
$$ p^{n+1}e_d\equiv p(\sum_{\ell =0}^{p^n-1} Tr_{K_{d-1}/\mathbb Q_p}(\zeta_{p^d}^{\ell })(1+T)^{\ell})
\pmod{\omega_n(T)}.$$
Now recall that $Tr_{K_{d-1}/\mathbb Q_p}(\zeta_{p^d}^{\ell })=0$ if $v_p(\ell )< d-1.$ Therefore:
$$p^{n+1}e_d \equiv p^d (\sum_{\ell =0}^{p^{n-d+1}-1}(1+T)^{\ell p^{d-1}}Tr_{K_0/\mathbb Q_p}(\zeta_p^{\ell
}))\pmod{\omega_n(T)}.$$
But:
$$\sum_{\ell =0, \ell \equiv 0 \pmod{p}}^{p^{n-d+1}-1}(1+T)^{\ell p^{d-1}}Tr_{K_0/\mathbb Q_p}(\zeta_p^{\ell
})= (p-1)\frac{\omega_n(T)}{\omega_{d}(T)},$$
and
$$\sum_{\ell =0, \ell \not \equiv 0 \pmod{p}}^{p^{n-d+1}-1}(1+T)^{\ell p^{d-1}}Tr_{K_0/\mathbb
Q_p}(\zeta_p^{\ell })=\frac{\omega_n(T)}{\omega_d(T)}-\frac{\omega_n(T)}{\omega_{d-1}(T)}.$$
Thus:
$$p^{n+1}e_d \equiv
p^{d+1}\frac{\omega_n(T)}{\omega_d(T)}-p^{d}\frac{\omega_n(T)}{\omega_{d-1}(T)}\pmod{\omega_n (T)}.$$
The Lemma follows. $\diamondsuit $\par
\section{Mirimanoff's power series}\par
${}$\par
Recall that Mirimanoff has introduced the following polynomials in $\mathbb F_p[T]:$
$$\forall j, \, 1\leq j\leq p-1,\, \varphi_j(T)=\sum_{a=1}^{p-1} a^{j-1}T^a.$$
These polynomials have many beautiful   properties and we refer the interested reader to \cite{RIB},
\cite{TER} and \cite{ANG}. 
In this section, we will introduce some power series that are related to these polynomials.\par
${}$\par
Let $L$ be a finite extension of $\mathbb Q_p.$ Let $\pi_L$ be a prime of $L.$ Let $\theta \in
\widehat{\Delta }.$ Let $a\in O_L$ such that $a(a-1) \not \equiv 0\pmod{\pi_L}.$ For all $n\geq 0,$ we set:
$$M_n(\theta , a)=\sum_{k=1, k\not \equiv 0\pmod{p}}^{p^{n+1}-1} \frac{a^k}{a^{p^{n+1}}-1} \theta
(k)\omega^{-1}(k) \gamma_n(k) \, \in O_L[\Gamma_n].$$
\newtheorem{Lemma3}{Lemma}[section]
\begin{Lemma3} \label{lemma3}
$$(M_n(\theta , a))_{n\geq 0}\in \lim_{\leftarrow}O_L[\Gamma_n].$$
\end{Lemma3}
\noindent{\sl Proof}  We must prove that:
$$\forall n\geq 0,\,  Res_{n+1,n}(M_{n+1}(\theta ,a))= M_n(\theta ,a).$$
Now observe that:
$$Res_{n+1,n}(M_{n+1}(\theta ,a))=\sum_{k=1, k\not \equiv 0\pmod{p}}^{p^{n+1}-1} \theta
(k)\omega^{-1}(k) \gamma_n(k) \frac{1}{a^{p^{n+2}}-1}(\sum_{\ell =0}^{p-1} a^{k+\ell p^{n+1}}).$$
The Lemma follows.$ \diamondsuit$\par
By the above Lemma, $(M_n(\theta , a))_{n\geq 0}$ corresponds to a power series $M(T,\theta , a)\in
\Lambda_L.$ If $\theta =\omega^j,$ $1\leq j \leq p-1,$ observe that:
$$M(0,\theta , a)\equiv \frac{\varphi_j(a)}{a^p -1}\pmod{\pi_L}.$$
Therefore we call $M(T,\theta ,a)$ the Mirimanoff's power series attached to $\theta $ and $a.$\par
\newtheorem{Lemma4}[Lemma3]{Lemma}
\begin{Lemma4} \label{Lemma4}
$$M(T,\theta ,a)=-\theta (-1)\omega^{-1}(-1) M(T, \theta ,a^{-1}).$$
\end{Lemma4}
\noindent{\sl Proof}\par
We have:
$$M_n(\theta ,a)=\sum_{k=1, k\not \equiv 0\pmod{p}}^{p^{n+1}-1} \frac{a^{p^{n+1}-k}}{a^{p^{n+1}}-1} \theta
(p^{n+1}-k)\omega^{-1}(p^{n+1}-k) \gamma_n(p^{n+1}-k).$$
Thus:
$$M_n(\theta ,a)=-\sum_{k=1, k\not \equiv 0\pmod{p}}^{p^{n+1}-1} \frac{(a^{-1})^k}{(a^{-1})^{p^{n+1}}-1}
\theta (p^{n+1}-k)\omega^{-1}(p^{n+1}-k) \gamma_n(p^{n+1}-k).$$
The Lemma follows. $\diamondsuit$\par
\newtheorem{Theorem2}[Lemma3]{Theorem}
\begin{Theorem2} \label{Theorem2}
${}$\par
$M'(T,\theta , a)\equiv 0\pmod{\pi_L}$ if and only if $a\equiv -1\pmod{\pi_L}$ and $\theta $ is odd.\par
\end{Theorem2}
\noindent{\sl Proof} By Lemma \ref{Lemma4}, $M(T,\theta ,-1)=0$ if $\theta $ is odd. Thus, if $a\equiv -1
\pmod{\pi_L}$ and if $\theta $ is odd, we have $M'(T,\theta , a)\equiv 0 \pmod{\pi_L}.$\par
The proof of this Theorem is based on Sinnott's proof that the Iwasawa $\mu$-invariant vanishes for cyclotomic
$\mathbb Z_p$-extensions of abelian number fields (\cite{SIN}) as exposed in Washington's book(\cite{WAS},
paragraph 16.2).\par
Now, let's suppose that $M'(T, \theta , a )\equiv 0\pmod{\pi_L}.$ We have:
$$\forall n\geq 0, \, M(T,\theta , a)\equiv\sum_{k=1,k\not \equiv 0\pmod{p}}^{p^{n+1}-1}\frac{a^k}{a^{p^{n+1}}-1}\theta
(k) \omega^{-1}(k) (1+T)^{i(k)}\pmod{\omega_n(T)},$$
where $i(k)={\rm Log}_p(k)/{\rm Log}_p(1+p).$ Thus, for all $n\geq 1,$ we have:
$$ (1+T) M'(T,\theta , a)\equiv\sum_{k=1,k\not \equiv 0\pmod{p}}^{p^{n+1}-1} i(k)\frac{a^k}{a^{p^{n+1}}-1}\theta
(k) \omega^{-1}(k) (1+T)^{i(k)}\pmod{(p^n ,\omega_n(T))}.$$
Therefore, for all $n\geq 1,$ we get:
$$\sum_{k=1,k\not \equiv 0\pmod{p}}^{p^{n+1}-1} i(k)\frac{a^k}{a^{p^{n+1}}-1}\theta
(k) \omega^{-1}(k) (1+T)^{i(k)}\equiv 0 \pmod{(\pi_L,\omega_n(T))}.$$
Recall that $i(k)\equiv i(k') \pmod{p^n}$ if and only if $\frac{<k>-1}{p}\equiv \frac{<k'>-1}{p}\pmod{p^n}.$ Therefore changing $i(k)$
to $\frac{<k>-1}{p}$ permutes exponents modulo $p^n$ and do not affect divisibility by $\pi_L.$ Thus:
$$\forall n\geq 1, \sum_{k=1,k\not \equiv 0\pmod{p}}^{p^{n+1}-1}\frac{<k>-1}{p} \frac{a^k}{a^{p^{n+1}}-1}\theta
(k) \omega^{-1}(k) (1+T)^{ \frac{<k>-1}{p}}\equiv 0 \pmod{(\pi_L,\omega_n(T))}.$$
Let $\alpha \in \mu_{p-1}.$ For $n\geq 1,$ set:
$$h_n^{\alpha }(t)\equiv \sum_{k=1,k \equiv \alpha \pmod{p}}^{p^{n+1}-1}\frac{\alpha^{-1}k-1}{p} \frac{a^k}{a^{p^{n+1}}-1}\theta
(k) \omega^{-1}(k) (1+T)^{ \frac{\alpha^{-1}k-1}{p}}\pmod{(p^n,\omega_n(T))}.$$
Note that:
$$h_{n+1}^{\alpha }(T)\equiv h_n^{\alpha }(T)\pmod{(p^n, \omega_n(T))}.$$
Now recall that:
$$\Lambda_L\simeq  \lim_{\leftarrow} \frac{\Lambda_L}{(p^n, \omega_n(T))}.$$
Therefore, there exists $h_{\alpha }(T)\in \Lambda_L$ such that:
$$\forall n\geq 1, h_{\alpha }(T) \equiv h_n^{\alpha }(T)\pmod{(p^n, \omega_n(T))}.$$
Thus, we have:
$$\sum_{\alpha \in \mu_{p-1}}h_{\alpha }(T) \equiv 0 \pmod{\pi_L}.$$
And also:
$$\sum_{\alpha \in \mu_{p-1}}(1+T)h_{\alpha }((1+T)^p-1) \equiv 0 \pmod{\pi_L}.$$
Now, note that:
$$(1+T)h_n^{\alpha }((1+T)^p -1)\equiv  f_n^{\alpha }((1+T)^{\alpha^{-1}}-1)\pmod{(p^n, \omega_n(T))},$$
where
$$ f_n^{\alpha }(T)\equiv \sum_{k=1,k \equiv \alpha
\pmod{p}}^{p^{n+1}-1}\frac{\alpha^{-1}k-1}{p}
\frac{a^k}{a^{p^{n+1}}-1}\theta (k) \omega^{-1}(k) (1+T)^{ k}\pmod{(p^n,\omega_n(T))}.$$
For $\alpha \in \mu_{p-1},$ let $s_0(\alpha )\in \{ 1,\cdots ,p-1\}$ such that $\alpha \equiv s_0(\alpha )\pmod{p}.$ We have:
$$f_n^{\alpha }(T)\equiv \theta (\alpha )\omega^{-1}(\alpha )\sum_{\ell =0}^{p^n -1} 
\frac{\alpha^{-1}(s_0(\alpha )+p\ell )-1}{p}\frac{a^{s_0(\alpha )+p\ell }}{a^{p^{n+1}}-1} (1+T)^{s_0(\alpha )+p\ell } \pmod{(p^n,
\omega_n (T))}.$$
But:
$$\sum_{\ell =0}^{p^n -1} 
\frac{a^{s_0(\alpha )+p\ell }}{a^{p^{n+1}}-1} (1+T)^{s_0(\alpha )+p\ell } \equiv
\frac{ a^{s_0(\alpha )} (1+T)^{s_0 (\alpha )}}{a^p (1+T)^p -1}\pmod{(p^n ,\omega_n (T))},$$
and
$$\sum_{\ell =0}^{p^n -1} 
\ell \frac{a^{s_0(\alpha )+p\ell }}{a^{p^{n+1}}-1} (1+T)^{s_0(\alpha )+p\ell } \equiv
-a^p (1+T)^p \frac{ a^{s_0(\alpha )} (1+T)^{s_0 (\alpha )}}{(a^p (1+T)^p -1)^2}\pmod{(p^n ,\omega_n (T))}.$$
Set:
$$f_{\alpha }(T)=\theta (\alpha )\omega^{-1}(\alpha ) \frac{ a^{s_0(\alpha )} (1+T)^{s_0 (\alpha )}}{a^p (1+T)^p -1}(
\frac{\alpha^{-1}s_0(\alpha )-1}{p}-\alpha^{-1}\frac{a^p (1+T)^p}{a^p (1+T)^p -1}).$$
Then:
$$\forall n\geq 1, f_n^{\alpha }(T)\equiv f_{\alpha }(T)\pmod{(p^n, \omega_n(T))}.$$
We have obtained:
$$\sum_{\alpha \in \mu_{p-1}}f_{\alpha }((1+T)^{\alpha^{-1}}-1)\equiv 0\pmod{\pi_L}.$$
Observe that for all $\alpha \in \mu_{p-1},$ $f_{\alpha }(T)\in \Lambda_L\cap L(T).$ Thus, by Sinnott's Lemma (\cite{WAS}, Lemma 16.9),
for all $\alpha \in \mu_{p-1}$ there exists $c_{\alpha }\in O_L$ such that:
$$f_{\alpha }(T)+f_{-\alpha }((1+T)^{-1}-1)\equiv c_{\alpha }\pmod{\pi_L}.$$
Observe that:
$$f_{\alpha }(T)\equiv \theta (\alpha )\omega^{-1}(\alpha )\frac{a^{s_0(\alpha )}}{a^p-1}(\frac{\alpha^{-1}s_0(\alpha
)-1}{p}-\alpha^{-1}\frac{a^p}{a^p-1})(1+T)^{s_0(\alpha )}\pmod{(\pi_L, T^p)},$$
and\par
\noindent $f_{-\alpha }((1+T)^{-1}-1)\equiv \theta (-\alpha )\omega^{-1}(-\alpha )\frac{a^{p-s_0(\alpha
)}}{a^p-1}(\frac{-\alpha^{-1}(p-s_0(\alpha ))-1}{p}+\alpha^{-1}\frac{a^p}{a^p-1})(1+T)^{s_0(\alpha )}\pmod{(\pi_L, T^p)}.$\par
\noindent Thus, for all $\alpha \in \mu_{p-1},$ we must have:\par
\noindent $a^{s_0(\alpha )}(\frac{\alpha^{-1}s_0(\alpha
)-1}{p}-\alpha^{-1}\frac{a^p}{a^p-1})\equiv \theta (-1)a^{p-s_0(\alpha )}(\frac{-\alpha^{-1}(p-s_0(\alpha
))-1}{p}+\alpha^{-1}\frac{a^p}{a^p-1})\pmod{\pi_L}.$\par
\noindent For $\alpha =1,$ we get:
$$a\frac{a^p}{a^p-1}\equiv -\theta (-1)a^{p-1}(\frac{a^p}{a^p-1}-1)\pmod {\pi_L}.$$
Thus:
$$a^2\equiv -\theta (-1)\pmod{\pi_L}.$$
We obtain $a\equiv -1\pmod{\pi_L}$ and $\theta $ is odd or $\theta $ is even and $a^2\equiv -1\pmod{\pi_L}.$ For the second case, if we
consider all the equations obtained when $\alpha $ runs through $\mu_{p-1},$ we obtain that for all $b\in \{ 1, \cdots ,p-1\},$ $b $
even, we must have:
$$b\equiv \omega (b)+p\pmod{p^2},$$
and for all $b\in \{ 1, \cdots ,p-1\},$ $b $
odd, we must have:
$$b\equiv \omega (b)\pmod{p^2},$$
This leads to a contradiction and the Theorem is proved. $\diamondsuit$\par
We will need the following Lemma:
\newtheorem{Lemma5}[Lemma3]{lemma}
\begin{Lemma5} \label{Lemma5}
There exists $\alpha \in \mu_{p-1}$ such that for all prime numbers  $\ell ,$ $\ell\equiv \alpha \pmod{p}$ and $\ell \geq p^2,$ we have:
$$Tr_{\mathbb Q(\zeta_{\ell})/\mathbb Q}(\frac{\zeta_{\ell}^{p+1}+\zeta_{\ell}^{p-1}}{(\zeta_{\ell}^p-1)^2})\not \equiv 0\pmod{p},$$
where $\zeta_{\ell }$ is a primitive $\ell$th root of unity.\par
\end{Lemma5}
\noindent{\sl Proof} For $a\in \mathbb Z,$ let $[a]_{\ell}\in \{ 0,\cdots ,\ell -1\}, $ such that $a\equiv [a]_{\ell }\pmod{\ell }.$
Let $\Phi_{\ell }(X)$ be the $\ell $th cyclotomic polynomial, then:
$$\Phi_{\ell }(1)=\ell ,$$
$$\Phi_{\ell }'(1)=\frac{\ell (\ell -1)}{2},$$
and
$$\Phi_{\ell}''(1)=\frac{\ell (\ell -1)(\ell -2)}{3}.$$
Now:
$$\frac{\Phi_{\ell}'(X)}{\Phi_{\ell}(X)}=\sum_{\rho\in \mu_{\ell}\setminus \{ 1\} }\frac{1}{X-\rho },$$
and
$$\frac{\Phi_{\ell}''(X)}{\Phi_{\ell}(X)}-(\frac{\Phi_{\ell}'(X)}{\Phi_{\ell}(X)})^2=-\sum_{\rho\in \mu_{\ell}\setminus \{ 1\}
}\frac{1}{(X-\rho )^2}.$$
Therefore:
$$Tr_{\mathbb Q(\zeta_{\ell })/\mathbb Q}(\frac{1}{\zeta_{\ell }-1})=\frac{1-\ell }{2},$$
and
$$Tr_{\mathbb Q(\zeta_{\ell })/\mathbb Q}(\frac{1}{(\zeta_{\ell }-1)^2})=\frac{(\ell -1)^2}{4}-\frac{(\ell-1)(\ell -2)}{3}.$$
Furthermore, if $a\in \mathbb Z,$ $a\not \equiv 0\pmod{\ell },$ we have:
$$Tr_{\mathbb Q(\zeta_{\ell })/\mathbb Q}(\frac{\zeta_{\ell}^a }{\zeta_{\ell }-1})=\frac{\ell +1}{2}-[a]_{\ell}.$$
Set:
$$S=Tr_{\mathbb Q(\zeta_{\ell})/\mathbb Q}(\frac{\zeta_{\ell}^{p+1}+\zeta_{\ell}^{p-1}}{(\zeta_{\ell}^p-1)^2}).$$
Let $m$ be the order of $p$ modulo $\ell,$   set:
$$a=[1+p^{m-1}]_{\ell}\, {\rm and }\, b=[1-p^{m-1}]_{\ell}.$$
We have:
$$S=Tr_{\mathbb Q(\zeta_{\ell})/\mathbb Q}(\frac{\zeta_{\ell}^{a}+\zeta_{\ell}^{b}}{(\zeta_{\ell}-1)^2}).$$
Since $\ell \geq p^2,$ we have $b\geq 3,$ and thus:
$$a=\ell +2-b.$$
Now:
$$S=Tr_{\mathbb Q(\zeta_{\ell})/\mathbb Q}(\frac{\zeta_{\ell}^a -1}{(\zeta_{\ell}^p-1)^2})+Tr_{\mathbb Q(\zeta_{\ell})/\mathbb
Q}(\frac{\zeta_{\ell}^b -1}{(\zeta_{\ell}^p-1)^2})+2Tr_{\mathbb Q(\zeta_{\ell})/\mathbb Q}(\frac{1}{(\zeta_{\ell}^p-1)^2}).$$
We obtain:
$$S=-b^2+b(\ell +2)-\frac{\ell^2+6\ell +5}{6}.$$
Therefore, $S\equiv 0\pmod{p}$ implies that:
$$\frac{\ell^2+2}{3}\in (\mathbb F_p)^2.$$
Now observe that the function field $\mathbb F_p(T,\sqrt{\frac{T^2+2}{3}})$ has genus zero and thus has exactly $p+1$ places of degree
one. The Lemma follows. $\diamondsuit$\par
\newtheorem{Theorem3}[Lemma3]{Theorem}
\begin{Theorem3} \label{Theorem3}
Let $\theta \in \widehat{\Delta },$ $\theta $ even. There exists $\alpha \in \mu_{p-1}$ such that for all prime numbers  $\ell ,$
$\ell\equiv \alpha \pmod{p}$ and $\ell \geq p^2,$ we have:
$$\sum_{\rho \in \mu_{\ell }\setminus\{ 1 \}}M'(T,\theta , \rho )\not \equiv 0 \pmod{p}.$$
\end{Theorem3}
\noindent{\sl Proof}\par
Let $\ell $ be a prime number, $\ell \not = p.$ Suppose that we have:
$$\sum_{\rho \in \mu_{\ell }\setminus\{ 1 \}}M'(T,\theta , \rho ) \equiv 0 \pmod{p}.$$
For $\alpha \in \mu_{p-1},$ set:
$$f_{\alpha }(T)=\sum_{\rho \in \mu_{\ell }\setminus\{ 1 \}}\theta (\alpha )\omega^{-1} (\alpha ) \frac{\rho^{s_0(\alpha )}
(1+T)^{s_0(\alpha )}}{\rho^p(1+T)^p-1}(\frac{\alpha^{-1}s_0(\alpha )-1}{p}-\alpha^{-1}\frac{\rho^p(1+T)^p}{\rho^p(1+T)^p-1}).$$
By the proof of Theorem \ref{Theorem2}, we get:
$$\sum_{\alpha \in \mu_{p-1}}f_{\alpha }((1+T)^{\alpha^{-1}}-1)\equiv 0\pmod{p}.$$
Therefore by \cite{WAS}, Lemma 16.9, for all $\alpha \in \mu_{p-1}$ there exists $c_{\alpha }\in \mathbb Z_p$ such that:
$$f_{\alpha }(T)+f_{-\alpha }((1+T)^{-1}-1)\equiv c_{\alpha }\pmod {p}.$$
But:
$$f_{1}(T)\equiv -(1+T)Tr_{\mathbb Q(\zeta_{\ell})/\mathbb Q}(\frac{\zeta_{\ell}^{p+1}}{(\zeta_{\ell}^p-1)^2})\pmod{(p, T^p)},$$
and
$$f_{-1}((1+T)^{-1}-1)\equiv -(1+T)Tr_{\mathbb Q(\zeta_{\ell})/\mathbb Q}(\frac{\zeta_{\ell}^{p-1}}{(\zeta_{\ell}^p-1)^2})\pmod{(p,
T^p)}.$$
Thus we get:
$$Tr_{\mathbb Q(\zeta_{\ell})/\mathbb Q}(\frac{\zeta_{\ell}^{p+1}+\zeta_{\ell}^{p-1}}{(\zeta_{\ell}^p-1)^2})\equiv 0\pmod{p}.$$
It remains to apply Lemma \ref{Lemma5}. $\diamondsuit$\par
\section{$p$-adic L-functions}
${}$\par
Let $\theta \in \widehat{\Delta }.$ We set:\par
\noindent - $f(T, \theta )=0$ if $\theta $ is odd,\par
\noindent - if $\theta $ is even, $f(T,\theta )$ is the Iwasawa power series associated to the $p$-adic L-function $L_p(s, \theta )$
(see \cite{WAS}, paragraph 7.2.).\par
\noindent Recall that if $\theta \not = 1,$ $f(T, \theta )\in \Lambda ,$ anf if $\theta $ is the trivial character then:
$$(1-\frac{1+p}{1+T})f(T,\theta )\in \Lambda^*.$$
For $\alpha \in \mu_{p-1},$ we set:
$$\rho_n(\alpha )=\frac{\alpha -\zeta_{p^{n+1}}}{\omega (\alpha -1)}\in U_n,$$
and 
$$\rho_{\infty}(\alpha )=(\rho_n(\alpha ))_{n\geq 0}\in U_{\infty}.$$
We set:
$$\eta_n=\prod_{\alpha\in \mu_{p-1}\setminus \{ 1 \}}\rho_n(\alpha )^{\alpha^{-1}-1}\in U_n,$$
and 
$$\eta_{\infty }=(\eta_n)_{n\geq 0}\in U_{\infty}.$$
\newtheorem{Lemma6}{Lemma}[section]
\begin{Lemma6} \label{Lemma6}
Let $\theta \in \widehat {\Delta },$ $\theta \not = 1, \omega.$ Then:
$$U_{\infty }^{e_{\theta }}=(\eta_{\infty }^{e_{\theta }})^{\Lambda }.$$
\end{Lemma6}
\noindent{\sl Proof} By the proof Theorem \ref{Theorem1}:
$$\phi_n^{\theta}(\rho_n(\alpha )^{e_{\theta }})\equiv M_n(\theta , \alpha^{-1})e_{\theta \omega^{-1}} T_n\pmod{p^{n+1}}.$$
Thus:
$$\phi_n^{\theta}(\eta_n^{e_{\theta}})\equiv (\sum_{\alpha \in \mu_{p-1}\setminus \{ 1 \}} (\alpha -1)M_n(\theta , \alpha ))e_{\theta
\omega^{-1}}T_n\pmod{p^{n+1}}.$$
But:
$$\sum_{\alpha \in \mu_{p-1}\setminus \{ 1 \}} (\alpha -1)M_0(\theta , \alpha )=(p-1)\theta (-1)\omega^{-1}(-1)\in \mathbb Z_p^*.$$
The Lemma follows. $\diamondsuit$\par
Let $\theta \in \widehat{\Delta },$ $\theta \not =1, \omega .$ We denote by $G(T,\theta )\in \Lambda^*$ the power series that
corresponds to $((\sum_{\alpha \in \mu_{p-1}\setminus \{ 1 \}} (\alpha -1)M_n(\theta , \alpha )))_{n\geq 0}.$ Observe that:
$$\sum_{\alpha \in \mu_{p-1}\setminus \{ 1 \}} (\alpha -1)M_n(\theta , \alpha )=(p-1)\theta (-1)\omega^{-1}(-1) \gamma_n(p-1)\sum_{\ell
=1, \ell \not \equiv 0 \pmod{p}}^{p^n}\theta (\ell )\omega^{-1}(\ell )\gamma_n(\ell ).$$
Let $\theta \in\widehat{\Delta },$ $\theta \not = 1,$ $\theta $ even. By $\pi_n^{e_{\theta}}$ we mean the unique $(p-1)^2$th root in
$U_n$ which is congruent to $1$ modulo $\pi_n$ of:
$$\prod_{a=1}^{p-1}(\frac{ \zeta_{p^{n+1}}^{\omega (a)}-1}{\zeta_{p^{n+1}}-1})^{(p-1) \theta^{-1}(a)}.$$
We set:
$$\pi_{\infty }^{e_{\theta }}=(\pi_n^{e_{\theta }})_{n\geq 0}\in U_{\infty}.$$
\newtheorem{Theorem4}[Lemma6]{Theorem}
\begin{Theorem4} \label{Theorem4}
${}$\par
\noindent i) Let $\theta \in \widehat{\Delta },$ $\theta $ even and non-trivial. Then:
$$(\pi_{\infty}^{e_{\theta }})^{-G(\frac{1+T}{1+p}-1,\theta )}=(\eta_{\infty}^{e_{\theta }})^{f(\frac{1+p}{1+T}-1,\theta )}.$$
ii)  Let $\theta \in \widehat {\Delta },$ $\theta \not =1,\omega .$ Then:
$$\forall \alpha \in \mu_{p-1}\setminus \{ 1 \} , (\rho_{\infty }(\alpha )^{e_{\theta }})^{G(\frac{1+T}{1+p}-1,\theta )}
=(\eta_{\infty}^{e_{\theta }})^{M(\frac{1+T}{1+p}-1, \theta , \alpha^{-1})}.$$
\end{Theorem4}
\noindent{\sl Proof} Recall that:
$$\phi_n^{\theta}(\rho_n(\alpha )^{e_{\theta }})\equiv M_n(\theta , \alpha^{-1})e_{\theta \omega^{-1}} T_n\pmod{p^{n+1}},$$
and
$$\phi_n^{\theta}(\eta_n^{e_{\theta}})\equiv (\sum_{\alpha \in \mu_{p-1}\setminus \{ 1 \}} (\alpha -1)M_n(\theta , \alpha ))e_{\theta
\omega^{-1}}T_n\pmod{p^{n+1}}.$$
Thus ii) follows from Lemma \ref{Lemma6}, Lemma \ref{Lemma2}, and Proposition \ref{Proposition1}.  Note that:
$$\phi_n^{\theta }(\pi_n^{e_{\theta }})\equiv e_{\theta \omega^{-1}} \frac{1}{\zeta_{p^{n+1}}-1}\pmod{p^{n+1}}.$$
Let $f(X)=\frac{X^{p^{n+1}}-1}{X-1},$ then:
$$(X-1)Xf'(X)+Xf(X)=p^{n+1}X^{p^{n+1}}.$$
Therefore:
$$\frac{1}{\zeta_{p^{n+1}}-1}=\frac{1}{p^{n+1}}\sum_{k=1}^{p^{n+1}-1}k\zeta_{p^{n+1}}^k.$$
Let $\theta \in \widehat {\Delta },$ for $\theta \not =1,$ set:
$$v_n(\theta )=-\frac{1}{p^{n+1}}\sum_{k=1,k\not \equiv 0 \pmod{p}}^{p^{n+1}-1}k\theta (k) \omega^{-1}(k) \gamma_n (k),$$
and if $\theta =1,$ set:
$$v_n(\theta )=-(1-(1+p)\gamma_n(1+p))\frac{1}{p^{n+1}}\sum_{k=1,k\not \equiv 0 \pmod{p}}^{p^{n+1}-1}k \omega^{-1}(k) \gamma_n (k).$$
Then by \cite{WAS}, paragraph 7.2, $(v_n(\theta ))_{n\geq 0}$ corresponds to $f(\frac{1}{1+T}-1,\theta )$ if $\theta \not =1,$ and to
$(1-(1+p)(1+T))f(\frac{1}{1+T}-1,\theta )$ if $\theta $ is trivial. Observe that if $\theta $ is even and non-trivial:
$$e_{\theta \omega^{-1}}\frac{1}{\zeta_{p^{n+1}}-1}=-v_n(\theta )e_{\theta \omega^{-1}} T_n,$$
and if $\theta =1 :$
$$(1-(1+p)\gamma_n(1+p))e_{\theta \omega^{-1}}\frac{1}{\zeta_{p^{n+1}}-1}=-v_n(\theta )e_{\theta \omega^{-1}} T_n.$$
The Theorem follows. $\diamondsuit$\par
\newtheorem{Theorem5}[Lemma6]{Theorem}
\begin{Theorem5} \label{Theorem5}
Let $d$ be an integer, $d\geq 2,$ $d\not \equiv 0 \pmod{p}.$ Let $\theta \in \widehat {\Delta }.$ Then:
$$\sum_{\rho \in \mu_d, \, o(\rho )=d}M(T, \theta , \rho) =
-f(\frac{1}{1+T}-1, \theta ) (\sum_{\ell \, \rm{divides} \, d}\ell \mu(\frac{d}{\ell })\theta (\ell ) \omega^{-1}(\ell )
(1+T)^{\frac{\rm{ Log}_p(\ell )}{\rm {Log}_p(1+p)}}),$$
wher $\mu (.)$ is the M\"obius function.\par
\end{Theorem5}
\noindent{\sl Proof} If $\theta $ is odd, the result is clear by Lemma \ref{Lemma4}. Thus we assume that $\theta $ is even. Let
$\Phi_d(X)$ be the $d$th cyclotomic polynomial, i.e.
$$\Phi_d(X)=\prod_{\rho \in \mu_d, \, o(\rho )=d}(X-\rho ) \in \mathbb Z[X].$$
Then:
$$\Phi_d(X)=\prod_{\ell \, \rm{divides} \, d}(X^{\ell }-1)^{\mu (d/\ell )}.$$
If we take logarithmic derivatives, we get:
$$\frac{\Phi_d'(X)}{\Phi_d(X)}=\sum_{\ell \, \rm{divides} \, d}\ell \mu( \frac{d}{\ell }) \frac{X^{\ell -1}}{X^{\ell }-1}.$$
It follows that:
$$\sum_{\rho \in \mu_d, \, o(\rho )=d} \frac{\zeta_{p^{n+1}}}{ \zeta_{p^{n+1}}-\rho } = \sum_{\ell \, \rm{divides} \, d} \ell \mu
(\frac{d}{\ell })\frac{\zeta_{p^{n+1}}^{\ell}}{\zeta_{p^{n+1}}^{\ell}-1}.$$
Thus:
$$\sum_{\rho \in \mu_d, \, o(\rho )=d} \frac{1}{ \rho \zeta_{ p^{n+1}}-1 } = (\sum_{\ell \, \rm{divides} \, d} \ell \mu
(\frac{d}{\ell })\sigma_n(\ell ))\frac{1}{\zeta_{p^{n+1}}-1}.$$
Now:
$$\frac{\rho^{p^{n+1}}-1}{\rho \zeta_{p^{n+1}}-1}-1=\sum_{k=1}^{p^{n+1}-1}\rho^k \zeta_{p^{n+1}}^k.$$
Thus:
$$\frac{1}{\rho \zeta_{p^{n+1}}-1}-\frac{1}{\rho^{p^{n+1}}-1}=\sum_{k=1}^{p^{n+1}-1}\frac{\rho^k }{\rho^{p^{n+1}}-1}\zeta_{p^{n+1}}^k.$$
We are working in the extension $\mathbb Q_p(\mu_d, \zeta_{p^{n+1}})/\mathbb Q_p(\mu_d )$  and we identify $\rm{Gal}(\mathbb Q_p(\mu_d,
\zeta_{p^{n+1}})/\mathbb Q_p(\mu_d ))$ with $\Gamma_n \times \Delta .$ Therefore:
$$e_{\theta \omega^{-1}}\frac{1}{\rho \zeta_{p^{n+1}}-1}=M_n(\theta , \rho) e_{\theta \omega^{-1}}\zeta_{p^{n+1}} + e_{\theta
\omega^{-1}}(\sum_{k=1}^{p^{n}-1}\frac{(\rho^p)^k }{(\rho^p)^{p^{n}}-1}\zeta_{p^{n}}^k).$$
Now observe that the map  $\mu_d\rightarrow \mu_d,$ $\rho \mapsto \rho^p,$ is an isomorphism. Thus:
$$e_{\theta \omega^{-1}}(\sum_{\rho \in \mu_d, \, o(\rho )=d} \frac{1}{\rho \zeta_{p^{n+1}}-1} )=(\sum_{\rho \in \mu_d, \, o(\rho
)=d}M_n(\theta , \rho ))e_{\theta \omega^{-1}}T_n.$$
But $T_n$ generates a normal basis for the field extension  $K_n/\mathbb Q_p,$ thus:
$$\sum_{\rho \in \mu_d, \, o(\rho )=d}M_n(\theta , \rho ) = -(\sum_{\ell \, \rm{divides} \, d} \ell \mu(\frac{d}{\ell }) \theta (\ell )
\omega^{-1}(\ell ) \gamma_n (\ell ))v_n (\theta ),$$
where $v_n (\theta )$ is as in the proof of Theorem \ref{Theorem4}. The Theorem follows. $\diamondsuit$\par
\newtheorem{Corollary1}[Lemma6]{Corollary}
\begin{Corollary1} \label{Corollary1}
Let $\theta \in \widehat {\Delta },$ $\theta$ even and non-trivial. Then:
$$f'(T,\theta )\not \equiv 0\pmod{p}.$$
\end{Corollary1}
\noindent{\sl Proof} Let $\alpha \in \mu_{p-1}$ as in Theorem \ref{Theorem3}. Let $\ell $ be  a prime number such that $\ell \geq p^2$
and $\ell \equiv \alpha \pmod{p^2}$ (note that there exist infinitely many such primes). We know that:
$$\sum_{\rho \in \mu_{\ell }\setminus\{ 1 \}}M'(T,\theta , \rho )\not \equiv 0 \pmod{p}.$$
But by Theorem \ref{Theorem5}:
$$\sum_{\rho \in \mu_{\ell }\setminus\{ 1 \}}M(T,\theta , \rho )=-f(\frac{1}{1+T}-1,\theta ) (\ell \theta (\ell )\omega^{-1} (\ell )
(1+T)^{\frac{\rm{ Log}_p(\ell )}{\rm {Log}_p(1+p)}}-1).$$
But since $\ell^{p-1}\equiv 1\pmod{p^2} ,$ we have:
$$\frac{\rm{ Log}_p(\ell )}{\rm {Log}_p(1+p)}\equiv 0\pmod{p}.$$
Thus, if we take derivatives and reduce modulo $p,$ we get:
$$\sum_{\rho \in \mu_{\ell }\setminus\{ 1 \}}M'(T,\theta , \rho ) \equiv \frac{1}{(1+T)^2}f'(\frac{1}{1+T}-1,\theta ) (\ell \theta
(\ell )\omega^{-1} (\ell ) (1+T)^{\frac{\rm{ Log}_p(\ell )}{\rm {Log}_p(1+p)}}-1)\pmod{p}.$$
The Corollary follows. $\diamondsuit$\par
The case of the trivial character is treated in the last section.
\section{Other results}
${}$\par
Let $\theta$ be an even Dirichlet character of conductor $d$ or $pd,$
where $d\geq 1,$ $d\not \equiv 0\pmod{p}.$ For all $n\geq 0,$ set:
$q_n=p^{n+1}d.$ Let $g(T,\theta )$ be the power series introduced in
\cite{WAS}, paragraph 7.2, i.e. 
$$g(T,\theta )=\frac{T-q_0}{1+T}f(T, \theta ),$$
where $f(T,\theta )$ is the power series associated to the $p$-adic
L-function $L_p(s,\theta )$ (see \cite{WAS}, Theorem 7.10). Set
$L=\mathbb Q_p (\theta ),$ and let $\pi_L$ be a prime of $L.$ Then, the
Ferrero-Washington Theorem states (see \cite{WAS}, paragraph 16.2):
$$g(T, \theta )\not \equiv 0 \pmod{\pi_L}.$$
We have:\par
\newtheorem{Theorem6}{Theorem}[section]
\begin{Theorem6} \label{Theorem6}
$$g'(T, \theta )\not \equiv 0 \pmod{\pi_L}.$$
\end{Theorem6}
\noindent{\sl Proof} For $y \in \mathbb Q,$ set:
$$B(y)=(1+q_0)\{ y\} -\{ (1+q_0)y\} -\frac{q_0}{2} \in \mathbb Z_p,$$
where $\{ y\}$ is the fractional part of $y.$ recall that (\cite {WAS},
proof of Theorem 7.10):
$$\forall n\geq 0,\,  g(T,\theta )\equiv 
\sum_{a=1, (a,q_0)=1}^{q_n} B(\frac{a}{q_n})\,  \theta \omega^{-1} (a)
\, (1+T)^{-i(a)-1}\pmod{\omega_n (T)},$$
where $i(a)=\frac{{\rm Log}_p (a)}{{\rm Log}_p(1+q_0)}.$ Let's suppose
that $g'(T, \theta )\equiv 0 \pmod{p}.$ We have, for all $n \geq 1:$\par
\noindent $ (1+T)^2g'(T,\theta )\equiv 
-\sum_{a=1, (a,q_0)=1}^{q_n}(i(a)+1)\,  B(\frac{a}{q_n})\,  \theta
\omega^{-1} (a)
\, (1+T)^{-i(a)}\pmod{(p^n, \omega_n (T))}.$\par
\noindent  Therefore, we get:
$$\forall n\geq 1, \sum_{a=1, (a,q_0)=1}^{q_n}(i(a)+1)\,  B(\frac{a}{q_n})\,  \theta
\omega^{-1} (a)
\, (1+T)^{i(a)}\equiv 0\pmod{(\pi_L, \omega_n (T))}.$$
Now, changing $i(a)$ to $(1+q_0)\frac{<a>-1}{p}$ permutes exponents
modulo $p^n$ and does not affect divisibility by $\pi_L.$ Thus, for all
$n\geq 1:$\par
\noindent $ \sum_{a=1,
(a,q_0)=1}^{q_n}((1+q_0)\frac{<a>-1}{p}+1)\,  B(\frac{a}{q_n})\, 
\theta
\omega^{-1} (a)
\, (1+T)^{(1+q_0)\frac{<a>-1}{p}}\equiv 0\pmod{(\pi_L, \omega_n (T))}.$\par
\noindent For $n\geq 1$ and for $\alpha \in \mu_{p-1},$ set:\par
\noindent $H_{\alpha }^n(T)\equiv \sum_{ a\equiv \alpha
\pmod{p}}((1+q_0)\frac{\alpha^{-1}a-1}{p}+1)\,  B(\frac{a}{q_n})\, 
\theta
\omega^{-1} (a)
\,  (1+T)^{(1+q_0)\frac{\alpha^{-1}a-1}{p}}\pmod{(p^n ,
\omega_n (T))}.$\par
\noindent Note that:
$$H_{\alpha }^{n+1}(T)\equiv H_{\alpha }^n (T)\pmod{(p^n, \omega_n
(T))}.$$
Thus, there exists $H_{\alpha }(T)\in \Lambda_L$ such that:
$$H_{\alpha }(T)\equiv H_{\alpha }^n (T)\pmod{(p^n, \omega_n
(T))}.$$
We get:
$$\sum_{\alpha \in \mu_{p-1}}H_{\alpha }(T) \equiv 0\pmod{\pi_L}.$$
Therefore:
$$\sum_{\alpha \in \mu_{p-1}}(1+T)^{1+q_0}H_{\alpha }((1+T)^p-1) \equiv
0\pmod{\pi_L}.$$
Let  $f_{\alpha }(T) \in \Lambda _L$ as in \cite{WAS}, Lemma 16.8. Set:
$$F_{\alpha }(T)= \frac{\alpha ^{-1}}{p} (1+T) f_{\alpha }' (T) +
(1-\frac{1+q_0}{p})f_{\alpha }(T).$$
It is not difficult to see that $F_{\alpha }(T) \in \Lambda_L \cap L(T).$
By \cite{WAS}, Lemma 16.8, we have:\par
\noindent $F_{\alpha }(T) \equiv \sum_{ a\equiv \alpha
\pmod{p}}((1+q_0)\frac{\alpha^{-1}a-1}{p}+1)\,  B(\frac{a}{q_n})\, 
\theta
\omega^{-1} (a)
\,  (1+T)^{(1+q_0)a}\pmod{(p^n ,
\omega_n (T))}.$\par
\noindent Therefore, we get:
$$\sum_{\alpha \in \mu_{p-1}}F_{\alpha }((1+T)^{\alpha^{-1}}-1) \equiv 0
\pmod{ \pi_L }.$$
Now apply \cite{WAS}, Lemma 16.9, and we get that for all $\alpha \in
\mu_{p-1},$ there exists $b_{\alpha } \in O_L$ such that:
$$F_{\alpha }(T)+F_{-\alpha }((1+T)^{-1}-1)\equiv b_{\alpha
}\pmod{\pi_L}.$$
But observe that:
$$F_{\alpha }(T) =F_{-\alpha }((1+T)^{-1}-1).$$
Thus, for all $\alpha \in \mu_{p-1}, $ there exists $c_{\alpha } \in O_L$
such that:
$$\frac{\alpha ^{-1}}{p} (1+T) f_{\alpha }' (T) +
(1-\frac{1+q_0}{p})f_{\alpha }(T)\equiv c_{\alpha }\pmod{\pi_L}.$$
Now we take $\alpha =1,$ and we set:\par
\noindent $G_1(T)=(1+q_0)\sum_{0<a<q_0, a\equiv 1\pmod{p}}\theta \omega^{-1}(a)\,
(1+T)^{a(1+q_0)}-\sum_{0<a<q_0^2+q_0, a\equiv 1\pmod{p}} \theta \omega^{-1} (a)\, (1+T)^a .$\par
\noindent We get:\par

\noindent $(1-\frac{1+q_0}{p})((1+T)^{q_0 (1+q_0)}-1)G_1(T)+\frac{1}{p}(1+T)
(((1+T)^{q_0(1+q_0)}-1)G_1'(T)-q_0(1+q_0)(1+T)^{q_0(1+q_0)-1}G_1(T))
 \equiv
c_1((1+T)^{q_0(1+q_0)}-1)^2\pmod{\pi_L}.$\par
Note that:
$$G_1(T)\equiv -(1+T) \pmod{(1+T)^{1+p}},$$
ans
$$G_1'(T)\equiv -1\pmod{(1+T)^p}.$$
Thus we must have:
$$c_1\equiv 0\pmod{\pi_L}.$$
The coefficient of $1+T$ in the left-side is: $1-d$ and the coefficient of $(1+T)^{q_0(1+q_0)+1}$ in the left side
is $2d-1+\frac{q_0^2}{p}.$ The Theorem follows. $\diamondsuit$\par

\end{document}